\documentclass[12pt]{article}
\scrollmode

\usepackage{amsfonts}
\usepackage{amsmath}
\usepackage{amssymb}
\usepackage{graphicx}

\allowdisplaybreaks

\def\eps{\varepsilon}
\def\qed{\hfill\rule{.2cm}{.2cm}}

\def\E{{\mathbb E}}
\def\Z{{\mathbb Z}}

\def\R{{\mathbb R}}
\def\C{{\mathbb C}}

\def\1{{\mathbf 1}}

\def\a{\alpha}
\def\o{\omega}  
\def\l{\lambda}

\newcommand{\AAA}          {\mathcal{A}}
\newcommand{\CC}          {\mathcal{C}}

\newcommand{\EE}          {\mathcal{E}}

\newcommand{\LL}         {\mathcal{L}}
\newcommand{\HH}         {\mathcal{H}}

\newcommand{\tX}        {{X^\xi}}
\newcommand{\tXn}        {{X^{\xi\,,\,n}}}
\newcommand{\tLLo}        {\LL^{\xi,\,\o}}

\newcommand{\tEEon}        {\EE^{\xi,\,\o,\,n}}
\newcommand{\oxi}        {\o^\xi}

\newcommand{\Var}       {\hbox{{\rm Var}}}

\def\nn{\nonumber}

\newtheorem{theo}{Theorem}[section]
\newtheorem{prop}[theo]{Proposition}
\newtheorem{lm}[theo]{Lemma}

\newtheorem{rmk}[theo]{Remark}

\def\beq{\begin{equation}}
\def\eeq{\end{equation}}
\newcommand{\bei}{\begin{itemize}}
\newcommand{\eei}{\end{itemize}}
\newcommand{\ben}{\begin{enumerate}}
\newcommand{\een}{\end{enumerate}}
\newcommand{\beqn}{\begin{eqnarray}}
\newcommand{\beqnn}{\begin{eqnarray*}}
\newcommand{\eeqn}{\end{eqnarray}}
\newcommand{\eeqnn}{\end{eqnarray*}}
\newcommand{\brm}{\begin{rmk}}
\newcommand{\erm}{\end{rmk}}

\title{Quenched invariance principles for random walks with random conductances. }
\author{
P.~Mathieu \footnote{ Universit\'e de Provence, CMI, 39 rue Joliot-Curie, 13013
Marseille, FRANCE.
pierre.mathieu@cmi.univ-mrs.fr}
}

\begin{document}
\maketitle

\begin{abstract}
We prove an almost sure invariance principle for a random walker
among i.i.d. conductances in $\Z^d$, $d\geq 2$. We assume
conductances are bounded from above but we do not require that they are
bounded from below.
\end{abstract}



\section{Introduction}

We consider continuous-time, nearest-neighbor random walks among random
(i.i.d.) conductances in $\Z^d$, $d\geq 2$ and prove that they satisfy
an almost sure invariance principle.

\subsection{Random walks and environments}
For $x, y \in \Z^d$, we write: $x \sim y$ if $x$ and $y$ are
neighbors in the grid $\Z^d$ and let $\E_d$ be the set of
non-oriented nearest-neighbor pairs $(x,y)$.\\
An {\it environment} is a function $\omega:\E_d\rightarrow
[0,+\infty[$. Since edges in $\E_d$ are not oriented, i.e. we
identified the edge $(x,y)$ with the reversed edge $(y,x)$, it is
implicit in the definition that environments are symmetric i.e.
$\o(x,y)=\o(y,x)$ for any pair of neighbors $x$ and $y$. \\
We let
$(\tau_z\,,\, z\in\Z^d)$ be the group of transformations of 
environments 
defined by $\tau_z\o(x,y)=\o(z+x,z+y)$.

We shall always assume that our environments are uniformly bounded from above.
Without loss of generality, we may assume that $\o(x,y)\leq 1$ for any edge.
Thus, for the rest of this paper, an environment will rather be a function
$\omega:\E_d\rightarrow [0,1]$.
We use the notation $\Omega=[0,1]^{E_d}$ for the set of environments (endowed with
the product topology and the corresponding Borel structure).
The value of an environment $\o$ at a given edge is called the {\it conductance}.

Let $\o\in\Omega$. We are interested in the behavior of the random
walk in the environment $\o$. We denote with $D(\R_+,\Z^d)$ the
space of c\`ad-l\`ag $\Z^d$-valued functions on $\R_+$ and let
$X(t)$, $t\in\R_+$, be the coordinate maps from $D(\R_+,\Z^d)$ to
$\Z^d$. The space $D(\R_+,\Z^d)$ is endowed with the Skorokhod
topology, see \cite{kn:Bill} or \cite{kn:JS}. For a given $\omega\in [0,1]^{\E_d}$ and for $x\in\Z^d$,
let $P^\o_x$ be the probability measure on $D(\R_+,\Z^d)$ under
which the coordinate process is the Markov chain starting at
$X(0)=x$ and with generator \beqn\label{int:gen} \LL^\o f(x)=\frac
1{n^\o(x)}\sum_{y\sim x} \o(x,y) (f(y)-f(x))\,, \eeqn where
$n^\o(x)=\sum_{y\sim x} \o(x,y)$. If $n^\o(x)=0$, let $\LL^\o
f(x)=0$ for any function $f$.

The behavior of $X(t)$ under $P^\o_x$ can be described as follows:
starting from point $x$, the random walker waits for an exponential
time of parameter $1$ and then chooses at random one of its
neighbors to jump to according to the probability law
$\o(x,.)/n^\o(x)$. This procedure is then iterated with independent
hopping times.

We have allowed environments to take the value $0$ and it is clear
from the definition of the random walk that $X$ will only travel
along edges with positive conductances. This remark motivates the
following definitions: call a {\it cluster} of the environment $\o$
a connected component of the graph $(\Z^d,\{e\in E_d\,;\,
\o(e)>0\})$. By construction, our random walker never leaves the
cluster of $\o$ it started from. Since edges are not oriented, the
measures with weights $n^\o(x)$ on the possibly different clusters
of $\o$ are reversible.

\subsection{Random environments} Let $Q$ be a product probability measure on $\Omega$.
In other words, we will now pick environments at random, in such a
way that the conductances of the different edges form a family of
independent identically distributed random variables. $Q$ is of
course invariant under the action of $\tau_z$ for any $z\in\Z^d$.

The random variables $(\1_{\o(e)>0}\,;\, e\in E_d)$ are independent Bernoulli
variables with common parameter $q=Q(\o(e)>0)$. Depending on the value of $q$,
a typical environment chosen w.r.t. $Q$ may or may not have infinite clusters.
More precisely, it is known from percolation theory that there is a critical
value $p_c$, that depends on the dimension $d$, such that for $q<p_c$,
$Q$.a.s. all clusters of $\o$ are finite and for $q>p_c$, $Q$.a.s. there is a
unique infinite cluster. In the first case the random walk is almost surely confined
to
a finite set and therefore does not satisfy the invariance principle (or satisfies
a degenerate version of it with vanishing asymptotic variance). We shall therefore
assume that the law $Q$ is {\it super-critical} i.e. that $$q=Q(\o(e)>0)>p_c\,.$$
Then the event `the origin belongs to the infinite cluster' has a non vanishing
$Q$ probability and we may define the conditional law:
\beqnn
Q_0(.)=Q(.\,\vert\, \hbox{$0$ belongs to the infinite cluster})\,.
\eeqnn

\subsection{Annealed results} Part of the analysis of the behavior of random
walks in random environments can be done using the
{\it point of view of the particle}: we consider the random walk $X$ started
at the origin and look at the random process describing the environment shifted
by the position of the random walker i.e. we let
$\o(t)=\tau_{X(t)}\o$. Thus $(\o(t)\,,\, t\in\R_+)$ is a random process taking
its values in $\Omega$.

Let us also introduce the measure \beqnn
{\tilde Q}_0(A)=\frac{\int_A n^\o(0)dQ_0(\o)}{\int n^\o(0)dQ_0(\o)}\,.
\eeqnn
Observe that ${\tilde Q}_0$ is obviously absolutely
continuous with respect to $Q_0$.

We list some of the properties of the process $\o(.)$ as proved in \cite{kn:DFGW}:
\begin{prop} \label{propDeMasi} (Lemmata 4.3 and 4.9 in \cite{kn:DFGW})\\
The random process $\o(t)$ is Markovian under $P^\o_0$.
The measure ${\tilde Q}_0$ is reversible, invariant and ergodic with respect to
$\o(t)$.
\end{prop}

Based on this proposition, the authors of \cite{kn:DFGW} could
deduce that the random walk $X(t)$ satisfies the invariance
principle {\it in the mean}. Let us define the so-called {\it
annealed} semi-direct product measure

\beqnn Q_0.P_x^\o[\,F(\o,X(.))\,]=\int P_x^\o[\,F(\o,X(.))\,]\,dQ_0(\o)\,.\eeqnn

\begin{theo} \label{theoDeMasi} (Annealed invariance principle, \cite{kn:DFGW})\\
Consider a random walk with i.i.d. super-critical conductances.
Under $Q_0.P_0^\o$, the process $(X^\eps(t)=\eps X(\frac
t{\eps^2}),t\in\R_+)$  converges in law to a non-degenerate Brownian
motion with covariance matrix $\sigma^2Id$ where $\sigma^2$ is
positive.
\end{theo}

It should be pointed out that the result of  \cite{kn:DFGW} is in fact
much more general. On one hand, \cite{kn:DFGW} deals with random
walks with unbounded jumps, under a mild second moment condition.
Besides, a similar annealed invariance principle is in fact proved
for any stationary law $Q$ rather than just product measures.

The positivity of $\sigma^2$ is not ensured by the general results
of \cite{kn:DFGW}) but it can be proved using comparison with the
Bernoulli case, see Remark \ref{rem:positivity}.

\subsection{The almost sure invariance principle}
The annealed invariance principle is not enough to give a completely
satisfactory description of the long time behavior of the random
walk. It is for instance clear that the annealed measure
$Q_0.P_0^\o$ retains all the symmetries of the grid. In particular
it is invariant under reflections through  hyperplanes passing
through the origin. This is  not true anymore for the law of the
random walk in a given environment. Still, one would expect
symmetries to be restored in the large scale, for a given
realization of $\o$.

Our main result is the following almost sure version of Theorem \ref{theoDeMasi}:

\begin{theo} \label{theorem1} (Quenched invariance principle)\\
Consider a random walk with i.i.d. super-critical conductances.
$Q_0$ almost surely, under $P^\o_0$, the
process $(X^\eps(t)=\eps X(\frac t{\eps^2}),t\in\R_+)$  converges in law
as $\eps$ tends to $0$ to a non-degenerate
Brownian motion with covariance matrix $\sigma^2Id$ where $\sigma^2$
is positive and
does not depend on $\o$.
\end{theo}

\subsection{The Bernoulli case and other cases} The main difficulty in proving Theorem \ref{theorem1}
is the lack of assumption on a lower bound for the values of the conductances.
Indeed, if one assumes that almost any environment is bounded from below by a
fixed constant i.e. there exists a $\delta>0$ such that $Q(\o(e)<\delta)=0$ 
then the conclusion of Theorem \ref{theorem1} was already proved in \cite{kn:SS}
using the classical `corrector approach' adapted from \cite{kn:Ko}.

Another special case recently solved is the Bernoulli case: let us
assume that only the values $0$ and $1$ are allowed for the
conductances i.e. $Q$ is a product of Bernoulli measures of
parameter $q$. Remember that we assume that we are in the
supercritical regime $q>p_c$. An environment can then be also
thought of as a (unweighted) random sub-graph of the grid and our
random walk is the simple symmetric random walk on the clusters of
the environment, i.e. jumps are performed according to the uniform
law on the neighbors of the current position in the graph $\o$.

In the Bernoulli case, quenched invariance principles have been
obtained by various authors in \cite{kn:BB}, \cite{kn:MP} and
\cite{kn:SS}. These three works develop different approaches to
handle the lack of a positive lower bound for the conductances. They
have in common the use of quantitative bounds on the transition
probabilities of the random walk. It is indeed known from
\cite{kn:Ba} that the kernel of the simple random walk on an
infinite percolation  cluster satisfies Gaussian bounds. A careful
analysis of the proofs shows that a necessary condition to obtain
the invariance principle using any of the three approaches in
\cite{kn:BB}, \cite{kn:MP} or \cite{kn:SS} is a Poincar\'e
inequality of the correct scaling (and in fact \cite{kn:MP} shows
that the Poincar\'e inequality is `almost' sufficient.) To be more
precise, let $A_n$ be the Poincar\'e constant on a box of size $n$
centered at the origin. In other words, $A_n$ is the inverse
spectral gap of the operator $\LL^\o$ restricted to the connected
component at the origin of the graph $\o\cap[-n,n]^d$  and with
reflection boundary conditions. Then one needs know that $Q_0$
almost surely, \beqn \label{int:poinc} \limsup n^{-2}A_n<\infty\,.
\eeqn Such a statement was originally proved in \cite{kn:MR} for the
Bernoulli case.

It turns out that (\ref{int:poinc}) is false in the general case of
i.i.d. conductances, even if one assumes that conductances are
always positive. We can choose for instance a product law with a
polynomial tail at the origin i.e. we assume  that there exists a
positive parameter $\gamma$ such that $Q(\o(e)\leq a)\sim a^\gamma$
as $a$ tends to $0$. Then it is not difficult to prove that, for
small values of $\gamma$, \beqnn \liminf\frac{\log A_n}{\log n}>
2\,.\eeqnn In \cite{kn:FM}, we considered a slightly different model
of symmetric random walks with random conductances with a polynomial
tail but non i.i.d. (although with finite range dependency only) and
we proved that \beqnn \frac{\log A_n}{\log n}\rightarrow 2\vee\frac
d\gamma\,,\eeqnn showing that, at least in the case $\gamma<d/2$,
the Poincar\'e constant is too big to be directly used to prove the
diffusive behavior of the random walk and one needs some new
ingredient to prove Theorem \ref{theorem1}.

\begin{rmk} In \cite{kn:FM}, we derived annealed estimates 
on the decay of the return probability of the random walk. 
More interestingly, in the very recent work \cite{kn:BBHK}, the authors 
could also obtain 
quenched  bounds on the decay of the return probability for quite general 
random walks with random conductances. 
Their results in particular show that anomalous decays do occur in high dimension. 
In such situations, although the almost sure invariance principle holds, see 
Theorem \ref{theorem1}, the local CLT fails.  
\end{rmk}

Our proof of Theorem \ref{theorem1} uses a time change argument that
we describe in the next part of the paper.
\vskip.5cm

{\it Acknowledgments:} the author would like to thank the referees of the first version 
of the paper for their careful reading and comments that lead to an improvement of the paper. 

{\it Note: after this paper was posted on the Arxiv, M. Biskup and 
T. Prescott wrote a preprint with a different proof of Theorem \ref{theorem1}, 
see \cite{kn:BiPres}. 
Their approach is based on ideas from \cite{kn:BB} when we prefer to invoke \cite{kn:MP}.   
They also need a time change argument, as here, and percolation results  
like Lemma \ref{lem:site''}.  
}


\vskip 1cm
\section{A time changed process}
\setcounter{equation}{0}
\label{sec:timechange}

In this section, we introduce a time changed process, $X^\xi$, and state an invariance
principle for it: Theorem \ref{theorem'}.

Choose a threshold parameter $\xi>0$ such that $Q(\o(e)\geq\xi)>p_c$.
For $Q$ almost any environment $\o$, the percolation graph
$(\Z^d,\{ e\in E_d\,;\, \o(e)\geq\xi\})$ has
a unique infinite cluster that we denote
with $\CC^\xi(\o)$.

By construction $\CC^\xi(\o)$ is a subset of $\CC(\o)$.
We will refer to the connected components of the complement of $\CC^\xi(\o)$ in
$\CC(\o)$ as {\it holes}.
By definition, holes are connected sub-graphs of the grid. Let $\HH^\xi(\o)$ be the
collection
of all holes. Note that holes may contain edges such that $\o(e)\geq\xi$.

We also define the conditioned measure
\beqnn Q_0^\xi(.)=Q(.\vert 0\in\CC^\xi(\o))\,.\eeqnn

Consider the following additive functional of the random walk:
\beqnn
A^\xi(t)=\int_0^t \1_{X(s)\in\CC^\xi(\o)}\,ds\,,
\eeqnn
its inverse $(A^\xi)^{-1}(t)=\inf \{s\,;\, A^\xi(s)>t\}$ and
define the corresponding time changed process
\beqnn
\tX(t)=X((A^\xi)^{-1}(t))\,.
\eeqnn

Thus the process $\tX$ is obtained by suppressing in the trajectory of $X$ all the
visits to the holes.
Note that, unlike $X$, the process $\tX$ may perform long jumps when straddling holes.

As $X$ performs the random walk in the environment $\o$, the
behavior of the random process $\tX$ is described in the next

\begin{prop}
Assume that the origin belongs to $\CC^\xi(\o)$. Then, 
under $P^\o_0$, the random process $\tX$ is a symmetric Markov
process on $\CC^\xi(\o)$.
\end{prop}

The Markov property, which is not difficult to prove, follows from a very general argument 
about time changed Markov processes. The reversibility of $\tX$ is a consequence of the 
reversibility of $X$ itself as will be discussed after equation (\ref{2:rates}). 

The generator of the process $\tX$ has the form
\beqn\label{2:gen'}
{\tLLo} f(x)=\frac 1{n^\o(x)}\sum_{y} \oxi(x,y) (f(y)-f(x))\,,
\eeqn
where
\beqn \nn
\frac{\oxi(x,y)}{n^\o(x)} &=&\lim_{t\rightarrow 0}  \frac 1t P_x^\o(\tX(t)=y)\\
&=&P_x^\o(\hbox{ $y$ is the next point in $\CC^\xi(\o)$ visited by the random walk
$X$})\,,
\label{2:rates} \eeqn
if both $x$ and $y$ belong to $\CC^\xi(\o)$ and $\oxi(x,y)=0$ otherwise.

The function $\oxi$ is symmetric: $\oxi(x,y)=\oxi(y,x)$ as follows
from the reversibility of $X$ and formula (\ref{2:rates}), but it is
no longer of nearest-neighbor type i.e. it might happen that
$\oxi(x,y)\not=0$ although $x$ and $y$ are not neighbors. More
precisely, one has the following picture: $\oxi(x,y)=0$ unless
either $x$ and $y$ are neighbors and $\o(x,y)\geq \xi$, or there exists a hole, $h$, such
that both $x$ and $y$ have  neighbors in $h$. (Both conditions may
be fulfilled by the same pair $(x,y)$.)

Consider a pair of neighboring points $x$ and $y$, both of them
belonging to the infinite cluster $\CC^\xi(\o)$ and such that
$\o(x,y)\geq\xi$, then \beqn \label{2:lowbo}
\oxi(x,y)\geq\xi\,.\eeqn This simple remark will play an important
role. It implies, in a sense to be made precise later, that the
parts of the trajectory of $\tX$ that consist in nearest-neighbors
jumps are similar to what the simple symmetric random walk on
$\CC^\xi(\o)$ does.

Finally observe that the environment $\oxi$ is stationary i.e. the law of $\oxi$
under $Q$
is invariant with respect to $\tau_z$ for all $z\in\Z^d$ as can be immediately seen
from
formula \ref{2:rates}.

\begin{theo} \label{theorem'} (Quenched invariance principle for $\tX$)\\
There exists a value $\xi_0>0$ such that for any $0<\xi\leq\xi_0$ the following holds.
For $Q_0$ almost any environment, under $P^\o_0$, the
process $(X^{\xi,\,\eps}(t)=\eps \tX(\frac t{\eps^2}),t\in\R_+)$  converges in law
as $\eps$ tends to $0$ to a non-degenerate
Brownian motion with covariance matrix $\sigma^2(\xi) Id$ where $\sigma^2(\xi)$
is positive and
does not depend on $\o$.
\end{theo}

The proof of Theorem \ref{theorem'} will be given in part
\ref{sec:proof}. It very closely mimics the arguments of
\cite{kn:MP}. Indeed, one uses the lower bound (\ref{2:lowbo}) to
bound the Dirichlet form of the process $\tX$ in terms of the
Dirichlet form of the simple symmetric random walk on $\CC^\xi(\o)$
and thus get the Poincar\'e inequality of the correct order. It is
then not difficult to adapt the approach of \cite{kn:MR} and
\cite{kn:Ba} to derive the tightness of the family $X^{\xi,\,\eps}$ 
and the  invariance principle follows as in
\cite{kn:MP}.

\begin{rmk} \label{rem:positivity}
The positivity of $\sigma^2$ in Theorem \ref{theorem1} and the
positivity of $\sigma^2(\xi)$ in Theorem \ref{theorem'} can be
checked using comparison arguments from \cite{kn:DFGW}. Indeed 
it follows from the expression of the effective diffusivity, 
see Theorem 4.5 part (iii) of \cite{kn:DFGW}, and from the 
discussion on monotonicity in part 3 of \cite{kn:DFGW} that $\sigma^2$
is an increasing function of the probability law $Q$ (up to some multiplicative factor). 
Therefore, if $Q$
stochastically dominates $Q'$ and the effective diffusivity under $Q'$ 
is positive, then the effective diffusivity under
$Q$ is also positive. 
Here $Q$ stochastically dominates the law of the environment with
conductances $\o'(e)=\xi\1_{\o(e)\geq\xi}$. The random walk in the
environment $\o'$ is the simple random walk on a percolation cluster
which is known to have a positive asymptotic diffusivity, see
\cite{kn:Ba} or the references in \cite{kn:MP}. The same argument
shows that $\sigma^2(\xi)>0$ for any $\xi$ such that
$Q(\o(e)\geq\xi)>p_c$.
\end{rmk}




To derive Theorem \ref{theorem1} from Theorem \ref{theorem'}, we
will compare the processes $X$ and $X^\xi$, for small values of
$\xi$. The large time asymptotic of the time change $A^\xi$ is easily
deduced from the ergodic theorem, as shown in Lemma \ref{lem:ergoA}
below and it implies that the asymptotic variance $\sigma^2(\xi)$ is
continuous at $\xi=0$, see Lemma \ref{lem:sigma}.

Let
\beqnn
c(\xi)={\tilde Q}_0 (0\in\CC^\xi(\o)) \,.
\eeqnn

\begin{lm} \label{lem:ergoA}
\beqnn
\frac {A^\xi(t)} t \rightarrow  c(\xi)\, \hbox{ $Q_0$ a.s.}
\eeqnn
as $t$ tends to $\infty$ and
 \beqn\label{2:cxi}
c(\xi)\rightarrow 1\,,
\eeqn
as $\xi$ tends to $0$.
\end{lm}
{\it Proof}: remember the notation $\o(t)=\tau_{X(t-)}\o$.
The additive functional $A^\xi(t)$ can also be written in the form
$A^\xi(t)=\int_0^t \1_{0\in\CC^\xi(\o(s))}\,ds$.

From Proposition \ref{propDeMasi}, we know that
${\tilde Q}_0$ is an invariant and ergodic measure for the process
$\o(t)=\tau_{X(t-)}\o$ and that it is absolutely continuous with respect to $Q_0$.

Thus the existence of the limit $\lim_{t\rightarrow +\infty} \frac
{A^\xi(t)} t$ follows from the ergodic theorem and the limit is
$c(\xi)={\tilde Q}_0(0\in\CC^\xi(\o))$. To check (\ref{2:cxi}), note
that $\1_{0\in\CC^\xi(\o)}$ almost surely converges to
$\1_{0\in\CC(\o)}$ as $\xi$ tends to $0$. Since ${\tilde Q}_0(0\in\CC(\o))=1$, we get that
$c(\xi)$ converges to $1$.

\qed

\begin{lm} \label{lem:sigma}
The asymptotic variances $\sigma^2$ in Theorem \ref{theoDeMasi} and $\sigma^2(\xi)$
from
Theorem \ref{theorem'}, and the constant $c(\xi)$ from Lemma \ref{lem:ergoA} satisfy
the
equality
\beqn\label{form:variance} c(\xi)\sigma^2(\xi)=\sigma^2\,.\eeqn
As a consequence, $\sigma^2(\xi)$ converges to $\sigma^2$ as $\xi$ tends to $0$.
\end{lm}

{\it Proof}:
formula (\ref{form:variance}) is deduced from Lemma \ref{lem:ergoA}. One can,
for instance, compute the law of the exit times from a large slab for both processes
$X$ and $X^\xi$.
Let $\tau(r)$ (resp. $\tau^\xi(r)$) be the exit time of $X$ (resp. $X^\xi$) from the
set
$[-r,r]\times\R^{d-1}$. Under the annealed measure, the Laplace transform of
$\tau(r)/r^2$
converges to $E(\exp(-\lambda T/\sigma^2))$ where $T$ is the exit time of $[-1,1]$
by a Brownian motion. This is a consequence of the invariance principle of Theorem
\ref{theoDeMasi}.
Theorem \ref{theorem'} implies that the Laplace transform of $\tau^\xi(r)/r^2$
converges
to $E(\exp(-\lambda T/\sigma^2(\xi)))$. (The convergence holds for $Q_0$ almost any
environment
and, by dominated convergence, under the annealed measure.) \\
On the other hand, we have $\tau^\xi(r)=A^\xi(\tau(r))$ and therefore Lemma
\ref{lem:ergoA}
implies that the Laplace transform of $\tau^\xi(r)/r^2$ has the same limit as the
Laplace
transform of $c(\xi)\tau^\xi(r)/r^2$ and therefore converges to $E(\exp(-\lambda
c(\xi)T/\sigma^2))$.
We deduce from these computations that
$$E(\exp(-\lambda c(\xi)T/\sigma^2))=E(\exp(-\lambda T/\sigma^2(\xi)))\,,$$
and, since this is true for any $\lambda\geq 0$, we must have
$c(\xi)\sigma^2(\xi)=\sigma^2$. \\
The continuity of
$\sigma^2(\xi)$ for $\xi=0$ is ensured by the continuity of $c(\xi)$.
\qed


\vskip 1cm
\section{How to deduce Theorem \ref{theorem1} from Theorem \ref{theorem'}}
\setcounter{equation}{0}
\label{sec:deduce}

We start stating a percolation lemma that will be useful to control
the contribution of holes to the behavior of the random walk.

\begin{lm} \label{lem:holes}
There exists a value $\xi_0>0$ such that for any $0<\xi\leq\xi_0$
the following holds. There exists a constant $a$ such that, $Q$
almost surely, for large enough $n$, the volume of any hole
$h\in\HH^\xi(\o)$ intersecting the box $[-n,n]^d$ is bounded from
above by $(\log n)^a$. ($a=7$ would do.)
\end{lm}

The proof of Lemma \ref{lem:holes} is postponed to part \ref{sec:perco}.

\subsection{Tightness}

In this section, we derive the tightness of the sequence of processes $X^\eps$ from
Theorem \ref{theorem'}.

\begin{lm} \label{lem:tight}
Under the assumptions of Theorem \ref{theorem1}, $Q_0$ almost
surely, under $P^\o_0$, the family of processes $(X^\eps(t)=\eps
X(\frac t{\eps^2}),t\in\R_+)$ is tight in the Skorokhod topology.
\end{lm}

{\it Proof}: we read from \cite{kn:JS}, paragraph 3.26, page 315 that a sequence of
processes
$x^\eps$ is tight if and only if the following two estimates hold:\\
(i) for any $T$, any $\delta>0$, there exist $\eps_0$ and $K$ such that for any
$\eps\leq\eps_0$
\beqn \label{ti1}
P(\sup_{t\leq T} \vert x^\eps(t)\vert\geq K)\leq\delta\,,
\eeqn
and\\
(ii) for any $T$, any $\delta>0$, any $\eta>0$, there exist $\eps_0$
and $\theta_0$ such that for any $\eps\leq\eps_0$ \beqn \label{ti2}
P(\sup_{v\leq u\leq T\,;\, u-v\leq\theta_0} \vert
x^\eps(u)-x^\eps(v)\vert>\eta)\leq\delta\,. \eeqn

Choose $\xi$ as in Theorem \ref{theorem'}. The sequence $X^{\xi,\,\eps}$ converges;
therefore it is tight and satisfies (\ref{ti1}) and (\ref{ti2}). By definition,
$$X^{\xi,\,\eps}(t)=X^\eps(\eps^2 (A^\xi)^{-1}(\frac t{\eps^2}))\,.$$

{\it Proof of condition (i)}: let us first check that $X^\eps$
satisfies (\ref{ti1}). \\ Assume that $\sup_{t\leq T} \vert
X^{\xi,\,\eps}(t)\vert\leq K$.  Given $t_0\leq T$, let
$x_0=X^\eps(t_0)$ i.e. $X(\frac {t_0}{\eps^2})=\frac 1 \eps {x_0}$
and define $s_0=\eps^2 A^\xi(\frac {t_0}{\eps^2})$. Since
$A^\xi(t)\leq t$, we have
$s_0\leq t_0$.\\
If $\frac 1 \eps {x_0}$ belongs to $\CC^\xi(\o)$, then $t_0=\eps^2
(A^\xi)^{-1}(\frac {s_0}{\eps^2})$
and $ X^{\xi,\,\eps}(s_0)=X^\eps(t_0)=x_0$ and therefore $\vert x_0\vert\leq K$. \\
Now suppose that $\frac 1 \eps {x_0}$ does not belong to $\CC^\xi(\o)$ and let
$t_1=\eps^2 (A^\xi)^{-1}(\frac {s_0}{\eps^2})$ and $x_1=X^\eps(t_1)$.
Then $t_1\leq t_0$ and $\frac 1 \eps {x_1}$ belongs to $\CC^\xi(\o)$. The same
argument as
before shows that $\vert x_1\vert\leq K$. On the other hand, by definition of the
time changed process $X^\xi$, $\frac 1\eps {x_1}$ is the last point in $\CC^\xi(\o)$
visited by $X$ before time $t_0$. Thus $\frac 1 \eps {x_0}$ belongs to a hole
on the boundary of which sits $\frac 1 \eps{x_1}$. It then follows from Lemma
\ref{lem:holes}
that $$\vert \frac 1 \eps{x_1}-\frac 1 \eps{x_0}\vert \leq (\log\frac K{\eps})^a\,.$$
Thus we have proved that
$$\vert x_0\vert\leq K+\eps (\log\frac K{\eps})^a\,.$$
We can choose $\eps_0$ small enough so that $\eps (\log\frac K{\eps})^a\leq K$
and therefore we have
$$\sup_{t\leq T} \vert X^{\xi,\,\eps}(t)\vert\leq K
\implies \sup_{t\leq T} \vert X^{\eps}(t)\vert\leq 2K\,.$$ Since the
sequence $X^{\xi,\,\eps}$ satisfies (\ref{ti1}), the event
`$\sup_{t\leq T} \vert X^{\xi,\,\eps}(t)\vert\leq K$' has a large
probability; therefore $\sup_{t\leq T} \vert X^{\eps}(t)\vert\leq
2K$ has a large probability and the sequence $X^{\eps}$ satisfies
(\ref{ti1}).

{\it Proof of condition (ii)}: as before, we will deduce that the sequence
$X^{\eps}$ satisfies (\ref{ti2}) from the fact that the sequence $X^{\xi,\,\eps}$
satisfies
(\ref{ti1}) and (\ref{ti2}).
Assume that
$$\sup_{v\leq u\leq T\,;\, u-v\leq\theta_0} \vert X^{\xi,\,\eps}(u)-
X^{\xi,\,\eps}(v)\vert\leq \eta\,.$$
We further assume that $\sup_{t\leq T} \vert X^{\xi,\,\eps}(t)\vert\leq K$.\\
Given $v_0\leq u_0\leq T$ such that $u_0-v_0\leq\theta_0$, let
$x_0=X^\eps(u_0)$, $y_0=X^\eps(v_0)$ and define $s_0=\eps^2
A^\xi(\frac {u_0}{\eps^2})$, $t_0=\eps^2 A^\xi(\frac
{v_0}{\eps^2})$, $u_1=\eps^2 (A^\xi)^{-1}(\frac {s_0}{\eps^2})$ and
$v_1=\eps^2 (A^\xi)^{-1}(\frac {t_0}{\eps^2})$.
Also let $x_1=X^\eps(u_1)$, $y_1=X^\eps(v_1)$.\\
Since $A^\xi(t)-A^\xi(s)\leq t-s$ whenever $s\leq t$, we have
$t_0\leq s_0\leq T$ and $s_0-t_0\leq\theta_0$. Besides, by
definition of $A^\xi$, we have $x_1= X^{\xi,\,\eps}(s_0)$ and $y_1=
X^{\xi,\,\eps}(t_0)$. We conclude that
$$ \vert x_1-y_1\vert\leq \eta\,.$$
On the other hand, the same argument as in the proof of condition (i) based on Lemma
\ref{lem:holes}
shows that
$$ \vert  x_1-x_0\vert+ \vert  y_1-y_0\vert\leq 2\eps  (\log\frac K{\eps})^a\,.$$
We have proved that
$$\sup_{v\leq u\leq T\,;\, u-v\leq\theta_0} \vert X^\eps(u)- X^\eps(v)\vert
\leq \eta+2\eps  (\log\frac K{\eps})^a\,.$$ Since both events
`$\sup_{v\leq u\leq T\,;\, u-v\leq\theta_0} \vert X^{\xi,\,\eps}(u)-
X^{\xi,\,\eps}(v)\vert\leq \eta$' and `$\sup_{t\leq T} \vert
X^{\xi,\,\eps}(t)\vert\leq K$' have large probabilities, we deduce
that  the processes $X^\eps$ satisfy condition (ii). \qed

\subsection{Convergence}

To conclude the derivation of Theorem \ref{theorem1} from Theorem \ref{theorem'}, it
only
remains to argue that, for any given time $t$, the two random variables
$X^\eps(t)$ and $X^{\xi,\,\eps}(t)$ are close to each other in probability.

\begin{lm} \label{lem:conv}
Under the assumptions of Theorem \ref{theorem1}, $Q_0$ almost surely,
for any $t$, any $\delta>0$, any $\eta>0$, then, for small enough $\xi$,
\beqnn
\limsup_{\eps\rightarrow 0} P^\o_0 ( \vert
X^\eps(t)-X^{\xi,\,\eps}(t)\vert>\eta)\leq\delta\,.
\eeqnn
\end{lm}

{\it Proof}: we shall rely on Lemma \ref{lem:ergoA}.
If $\vert X^\eps(t)-X^{\xi,\,\eps}(t)\vert>\eta$, then one of the following two
events must hold:
$$(I)=
\{\sup_{\theta c(\xi)t\leq s\leq t} \vert
X^{\xi,\,\eps}(s)-X^{\xi,\,\eps}(t)\vert>\frac \eta 2\}\,,$$
$$(II)
=\{ \inf_{\theta c(\xi)t\leq s\leq t} \vert X^{\xi,\,\eps}(s)-X^\eps(t)\vert>\frac
\eta 2\}\,.$$
Here $\theta$ is a parameter in $]0,1[$.\\
The invariance principle for $X^{\xi,\,\eps}$, see Theorem
\ref{theorem'}, implies that the probability of $(I)$ converges as
$\eps$ tends to $0$ to the probability $P(\sup_{\theta c(\xi)t\leq
s\leq t} \sigma(\xi)\vert B(s)-B(t)\vert>\frac \eta 2)$, where $B$
is a Brownian motion. Since $\sigma(\xi)$ is bounded away from $0$,
see Lemma \ref{lem:sigma}, and since $c(\xi)\rightarrow 1$ as
$\xi\rightarrow 0$, we deduce that there exists a value for $\theta$
such that \beqn\label{eq:12} \limsup_{\xi \rightarrow 0}\limsup_{\eps\rightarrow 0}
P^\o_0 (I)\leq \delta\,. \eeqn

We now assume that $\theta$ has been chosen so that (\ref{eq:12})
holds. We shall end the proof of the Lemma by showing  that \beqn
\label{eq:11} \limsup_{\eps\rightarrow 0} P^\o_0 (II)=0\,. \eeqn
Since, from the tightness of the processes $X^\eps$, see Lemma
\ref{lem:tight}, we have
$$\limsup_{\eps\rightarrow 0}P^\o_0 (\sup_{s\leq t}\vert X^\eps(s)\vert\geq
\eps^{-1})=0\,,$$
we will estimate the probability that both events $(II)$ and
`$\sup_{s\leq t}\vert X^\eps(s)\vert\leq \eps^{-1}$' hold. \\
Let $u=\eps^2 A^\xi(\frac t{\eps^2})$ and note that $u\leq t$. From Lemma
\ref{lem:ergoA},
we know that $u\geq \theta c(\xi) t$ for small enough $\eps$ depending on $\o$. \\
If $X^\eps(t)$ belongs to $\CC^\xi(\o)$, then $X^\eps(t)=X^{\xi,\,\eps}(u)$
and therefore $(II)$ does not hold.  \\
Otherwise $X^\eps(t)$ belongs to a hole on the boundary of which sits
$X^{\xi,\,\eps}(u)$.
Using the condition $\sup_{s\leq t}\vert X^\eps(s)\vert\leq \eps^{-1}$ and
Lemma \ref{lem:holes}, we get that
$$\vert X^\eps(t)-X^{\xi,\,\eps}(u)\vert\leq \eps(\log \frac 1\eps)^a\,.$$
For sufficiently small $\eps$ we have $\eps(\log 1/\eps)^a<\frac \eta 2$ and
therefore $(II)$ fails. The proof of (\ref{eq:11}) is complete.
\qed

{\it End of the proof of Theorem \ref{theorem1}}: choose times
$0<t_1<...<t_k$. Use Lemma \ref{lem:conv}, to deduce that for small
enough $\xi$, as $\eps$ tends to $0$, the law of
$(X^\eps(t_1),...,X^\eps(t_k))$ comes close to the law of
$(X^{\xi,\,\eps}(t_1),...,X^{\xi,\,\eps}(t_k))$, which in turn,
according to Theorem \ref{theorem'}, converges to the law of
$(\sigma(\xi)B(t_1),...,\sigma(\xi)B(t_k))$, where $B$ is a Brownian
motion. We now let $\xi$ tend to $0$: since $\sigma(\xi)$ converges
to $\sigma$, see Lemma \ref{lem:sigma}, the limiting law of
$(X^\eps(t_1),...,X^\eps(t_k))$ is the law of $(\sigma
B(t_1),...,\sigma B(t_k))$ i.e. we have proved that $X^\eps$
converges in law to a Brownian motion with variance $\sigma^2$ in
the sense of finite dimensional marginals. The tightness Lemma
\ref{lem:tight} implies that the convergence in fact holds in the
Skorokhod topology. \qed




\vskip 1cm
\section{Proof of Theorem \ref{theorem'}}
\setcounter{equation}{0}
\label{sec:proof}

We will outline here a proof of Theorem \ref{theorem'}.
Our strategy is quite similar to the one recently used in
\cite{kn:MR}, \cite{kn:Ba} and \cite{kn:MP} to study the simple symmetric
random walk on a percolation cluster. No new idea is required.

\vskip .5cm {\bf Step 0: notation}

As before, we use the notation $\omega$ to denote a typical environment
under the measure $Q$.
For a given edge $e\in\E_d$ (and a given choice of $\omega$), we define
\beqnn
\alpha(e)=\1_{\o(e)>0}\,;\, \alpha'(e)=\1_{\o(e)\geq\xi}\,.
\eeqnn
As in part \ref{sec:timechange}, let $\CC^\xi(\o)$ be the infinite cluster of the
percolation
graph $\alpha'$. For $x,y\in\CC^\xi(\o)$,
we define the {\it chemical distance} $d^\xi_\o(x,y)$ as the minimal number of jumps
required
for the process $\tX$ to go from $x$ to $y$, see part \ref{sec:deviation}.

We recall the definition of the generator ${\tLLo}$ from formula (\ref{2:gen'}).
Since the function $\oxi$ is symmetric, the operator ${\tLLo}$ is reversible with
respect
to the measure $\mu_\o=\sum_{z\in\CC^\xi(\o)} n^\o(z)\delta_z$.

Let  $\CC^{n}(\o)$ be the connected
component of $\CC^\xi(\o)\cap [-n,n]^d$ that contains the
origin. Let $(\tXn(t),\ t \geq 0)$ be the random walk $\tX$ restricted to the set
$\CC^n(\o)$.
The definition of $\tXn$ is the same as for $\tX$ except that jumps outside $\CC^n$
are
now forbidden. Its Dirichlet form is
\beqnn
\tEEon(f,f)=\frac 12 \sum_{x\sim y\in \CC^n(\o)}
                        \oxi(x,y) (f(x)-f(y))^2
\eeqnn

We use the notation $\tau^n$ for the exit time of the process $\tX$ from the box
$[-2n+1,2n-1]^d$ i.e. $\tau^n=\inf\{t\,;\, \tX(t)\notin[-2n+1,2n-1]^d\}$.

\vskip .5cm {\bf Step 1: Carne-Varopoulos bound}

The measure $\mu_\o$ being reversible for the process $\tX$,
the transition probabilities satisfy a Carne-Varopoulos bound:
\beqnn
P^\o_x(\tX(t)=y)\leq Ce^{-d^\xi_\o(x,y)^2/(4t)}+e^{-ct}\,,
\eeqnn
where $c=\log 4-1$ and $C$ is some constant that depends on $\xi$ and $\o$. (See \cite{kn:MR}, appendix C.)

By Lemma \ref{lem:deviation}, we  can replace the chemical distance
$d^\xi_\o(x,y)$ by the Euclidean distance $\vert x-y\vert$, provided
that $x\in[-n,n]^d$ and $n$ is large enough. We get that, $Q_0^\xi$
almost surely, for large enough $n$, for any $x\in[-n,n]^d$ and any
$y\in\Z^d$ such that $\vert x-y\vert\geq (\log n)^2$, then
\beqn\label{eq:carnevaro}  P^\o_x(\tX(t)=y)\leq 
Ce^{-\frac{\vert
x-y \vert^2}{Ct}}+e^{-ct}\,.\eeqn 

The same reasoning as in \cite{kn:MR}, appendix C (using Lemma
\ref{lem:deviation} again) then leads to upper bounds 
for the exit time $\tau^n$: $Q_0^\xi$
almost surely, for large enough $n$, for any $x\in[-n,n]^d$ and any
$t$, we have 
\beqn \label{eq:taun}
P_x^\omega[\tau^n\leq t] \leq Ct
n^d e^{-\frac{n^2}{Ct}}+e^{-ct}\,. \eeqn
Indeed, let $N(t)$ be the number of jumps the random walk performs until time $t$ 
and let $\sigma^n$ be the number of jumps of the walk until it exits the box $[-2n+1,2n-1]^d$, 
so that $\sigma^n=N(\tau^n)$. 
Note that the process $(N(t)\,,\, t\in\R_+)$ is a Poisson process of rate $1$. 
With probability larger than $1-e^{-ct}$, we have $N(t)\leq 2t$. 
If $N(t)\leq 2t$ and $\tau^n\leq t$, then $\sigma^n\leq 2t$ and there are  
at most $2t$ choices for the value of $\sigma^n$. 
Let $y$ be the position of the walk at the exit time and let $z$ be the last point  
visited before exiting.  Note that $d^\xi_\o(z,y)=1$. 
Due to Lemma
\ref{lem:deviation}, we have 
\beqnn 
\vert x-y\vert\leq \frac 1{c^-} d^\xi_\o(x,y) 
\leq \frac 1{c^-} (d^\xi_\o(x,z)+1) \leq \frac{c^+}{c^-} (1+\vert x-z\vert)\leq 
\frac{c^+}{c^-} (1+2n)\,.   
\eeqnn 
Note that our use of Lemma
\ref{lem:deviation} here is legitimate. Indeed $\vert x-y\vert$ is of order $n$ 
and, since $d^\xi_\o(z,y)=1$, Lemma \ref{lem:holes} implies that $\vert y-z\vert$ is at most 
of order $(\log n)^7$. Therefore $\vert x-z\vert$ is of order $n$ and thus certainly larger that $(\log n)^2$. 

Thus we see that there are at most of order $n^d$ possible choices for $y$.  
Finally, due to (\ref{eq:carnevaro}), 
\beqnn P^\o_x(\tX(s)=y)\leq C e^{-\frac{n^2}{Ct}}\,,\eeqnn 
for any $s\leq 2t$, $x\in[-n,n]^d$ and $y\notin [-2n+1,2n-1]^d$. 
Putting everything together, we get (\ref{eq:taun}). 

\vskip .5cm {\bf Step 2: Nash inequalities and on-diagonal decay}

\begin{lm}\label{lem:1}
For any $\theta>0$,  there exists a constant $c_u(\theta)$  such
that, $Q_0^\xi$ a.s. for large enough $t$, we have \beqn
\label{eq:ondiag}
 P_x^{\omega}[\tX(t)=y]
 \leq
 \frac{c_u(\theta)}{t^{d/2}}\,,
\eeqn for any $x\in \CC^\xi(\o)$ and $y\in\Z^d$ such that $\vert x\vert\leq
t^\theta$.
\end{lm}

{\it Proof}:

We use the notation $\alpha'(e)=\1_{\o(e)\geq\xi}$. Note that the random variables 
$(\alpha'(e)\,;\, e\in E_d)$ are independent Bernoulli
variables with common parameter $Q(\alpha'(e)>0)=Q(\o(e)\geq\xi)$. Since we have assumed that 
$Q(\o(e)\geq\xi)>p_c$, the environment $\alpha'$ is a typical realization of super-critical bond 
percolation.

The following Nash inequality is proved in \cite{kn:MR}, equation (5):
there exists a constant $\beta$ such that
$Q_0^\xi$ a.s. for large enough $n$, for any function
$f:\CC^n(\o)\rightarrow\R$ one has
\beqnn
 \Var(f)^{1+\frac{2}{\eps(n)}}
  \leq
\beta \,
  n^{2(1-\frac{d}{\eps(n)})}\,\EE^{\alpha',n}(f,\,f)\,\|f\|_1^{4/\eps(n)}\,,
\eeqnn where \beqnn \EE^{\alpha',n}(f,f)=\frac 12 \sum_{x\sim y\in
\CC^n(\o)}
                        \alpha'(x,y) (f(x)-f(y))^2
\,.
\eeqnn

The variance and the $L_1$ norms are computed with respect to the
counting measure on $\CC^n(\o)$ and $\eps(n)=d+2d\frac{\log\log
n}{\log n}$.
(Note that there is a typo in \cite{kn:MR} where it is claimed that 
(5) holds for the uniform probability on $\CC^n(\o)$ instead of the counting measure.) 

Inequality (\ref{2:lowbo}) implies that $\alpha'(x,y)\leq \xi^{-1}
\oxi(x,y)$. Therefore $\EE^{\alpha',n}$ and $\tEEon$ satisfy the
inequality \beqn\label{eq:compDF} \EE^{\alpha',n}(f,f)\leq \frac
1\xi\, \tEEon(f,f)\,. \eeqn

Using inequality (\ref{eq:compDF}) in the previous Nash inequality, we
deduce that there exists a constant $\beta$ (that depends on $\xi$) such that
$Q_0^\xi$ a.s. for large enough $n$, for any function
$f:\CC^n(\o)\rightarrow\R$ one has
\beqn \label{eq:nash}
 \Var(f)^{1+\frac{2}{\eps(n)}}
  \leq
\beta\,
  n^{2(1-\frac{d}{\eps(n)})}\,\tEEon(f,\,f)\,\|f\|_1^{4/\eps(n)}\,.
\eeqn

As shown in \cite{kn:MR} part 4, the Carne-Varopoulos inequality
(\ref{eq:carnevaro}), inequality (\ref{eq:taun}) and the Nash
inequality (\ref{eq:nash}) can be combined to prove upper bounds on
the transition probabilities. We thus obtain that: there exists a
constant $c_u$ such that, $Q_0^\xi$ a.s. for large enough $t$, we
have \beqn \label{eq:ondiag0}
 P_0^{\omega}[\tX(t)=y]
 \leq
 \frac{c_u}{t^{d/2}}\,,
\eeqn for any  $y\in\Z^d$.

Using the translation invariance of $Q$, it is clear that estimate
(\ref{eq:ondiag0}) in fact holds if we choose another point
$x\in\Z^d$ to play the role of the origin. Thus, for any $x\in\Z^d$,
$Q$ a.s. on the set $x\in \CC^\xi(\o)$, for  $t$ larger than some
random value $t_0(x)$, we have \beqn \label{eq:ondiag'}
 P_x^{\omega}[\tX(t)=y]
 \leq
 \frac{c_u}{t^{d/2}}\,,
\eeqn for any  $y\in\Z^d$.

In order to deduce the Lemma from the upper bound
(\ref{eq:ondiag'}), one needs control the tail of the law of
$t_0(0)$. \\
Looking at the proofs in \cite{kn:MR}, one sees that all the error
probabilities decay faster than any polynomial. More precisely, the
$Q_0^\xi$ probability that inequality (\ref{eq:nash}) fails for some
$n\geq n_0$ decays faster than any polynomial in $n_0$. From the
proof of Lemma \ref{lem:deviation}, we also know that the $Q_0^\xi$
probability that inequality (\ref{eq:carnevaro}) fails for some
$n\geq n_0$ decays faster than any polynomial in $n_0$. As a
consequence, a similar bound holds for inequality (\ref{eq:taun}).\\
To deduce error bounds for (\ref{eq:ondiag0}), one then needs to go
to part 4 of \cite{kn:MR}. Since the proof of the upper bound
(\ref{eq:ondiag0}) is deduced from (\ref{eq:carnevaro}),
(\ref{eq:taun}) and (\ref{eq:nash}) by choosing $t\log t =b n^2$ for
an appropriate constant $b$, we get that  $Q_0^\xi(\hbox {inequality
(\ref{eq:ondiag0}) fails for some $t\geq t_0$})$ decays faster than
any polynomial in $t_0$. By translation invariance, the same holds
for (\ref{eq:ondiag'}) i.e. for any $A>0$, there exists $T$ such
that 
\beqnn Q(x\in \CC^\xi(\o)\,\hbox{and}\,t_0(x)\geq t_0 )\leq
t_0^{-A}\,,\eeqnn 
for any $t_0>T$. Therefore, 
\beqnn Q(\exists x\in \CC^\xi(\o)\,;\,
\vert x\vert\leq t_0^\theta \,\hbox{and}\, t_0(x)\geq t_0)\leq t_0^{d\theta-A}\,.\eeqnn 
One then chooses $A$
larger than $d\theta+1$ and the Borel-Cantelli lemma gives the end of
the proof of (\ref{eq:ondiag}).\qed

\vskip .5cm {\bf Step 3: exit times estimates and tightness}

We denote with $\tau(x,r)$ the exit time of the random walk from the
ball of center $x$ and Euclidean radius $r$.

\begin{lm}\label{lem:2}
For any $\theta>0$,  there exists a constant $c_e$  such that,
$Q_0^\xi$ a.s. for large enough $t$, we have \beqn \label{eq:exit}
 P_x^{\omega}[\tau(x,r)< t]
 \leq
c_e \frac{\sqrt{t}}r\,, \eeqn for any $x\in \Z^d$ and $r$ such that
$\vert x\vert\leq t^\theta$ and $r\leq t^\theta$.
\end{lm}

{\it Proof}: the argument is the same as in \cite{kn:Ba}, part 3. We
define
$$ M_x(t)=E^\o_x[d^\xi_\o(x,X^\xi(t))]$$
and
$$ Q_x(t)=-E^\o_x[\log q^\o_t(x,X^\xi(t))]\,,$$
where $q^\o_t(x,y)=P^\o_x(X^\xi(t)=y)/\mu_\o(x)$. Then, for large
enough $t$ and for $\vert x\vert\leq t^\theta$, one has:
\beqnn &&Q_x(t)\geq -\log c_u+\frac d 2 \log t\,,\\
&&M_x(t)\geq c_2 \exp(Q_x(t)/d)\,,\\
&&Q_x'(t)\geq \frac 12 (M_x'(t))^2\,. \eeqnn The first inequality is
obtained as an immediate consequence of Lemma \ref{lem:1}. The
second one is proved as in \cite{kn:Ba}, Lemma 3.3 and the third one
as in \cite{kn:Ba}, equation (3.10), using ideas from \cite{kn:Bass} 
and \cite{kn:Nash}. 
Note that, in the proof of the
second inequality, we used Lemma \ref{lem:deviation} to control the
volume growth in the chemical distance $d^\xi_\o$. One now
integrates these inequalities to deduce that \beqn\label{eq:mean}
c_1\sqrt t \leq M_x(t)\leq c_2 \sqrt t\,. \eeqn Once again the proof
is the same as in \cite{kn:Ba}, Proposition 3.4. Note that, in the
notation of \cite{kn:Ba}, $T_B=\vert x\vert^{1/\theta}$ so that
equation (\ref{eq:mean}) holds for $t\geq \frac 1\theta \vert
x\vert^{1/\theta}\log\vert x\vert$. The end of the proof is
identical to the proof of Equation (3.13) in \cite{kn:Ba}. \qed

\begin{lm}\label{lem:2'}
$Q_0^\xi$ a.s. for large enough $t$, we have \beqn \label{eq:exit'}
 P_x^{\omega}[\tau(x,r)< t]
 \leq
27 
(c_e)^3 (\frac{\sqrt{t}}r)^3\,, \eeqn for any $x\in \Z^d$ and $r$ such that
$\vert x\vert\leq t^\theta$ and $r\leq t^\theta$.
\end{lm}

{\it Proof}: let $x'=X^\xi(\tau(x,r/3))$, $x''=X^\xi(\tau'(x',r/3))$ where $\tau'(x',r/3)$ is the exit 
time from the ball of center $x'$ and radius $r/3$ after time $\tau(x,r/3)$ and let $\tau''(x'',r/3)$ 
be the exit 
time from the ball of center $x''$ and radius $r/3$ after time $\tau'(x,r/3)$. 
In order that $\tau(x,r)< t$ under $P_x^{\omega}$ we must have 
$\tau(x,r/3)< t$  and 
$\tau'(x',r/3)< t$ and $\tau''(x'',r/3)<t$. We can then use Lemma \ref{lem:2} to estimate the probabilities 
of these $3$ events and conclude that (\ref{eq:exit'}) holds.  \qed

\begin{lm}\label{lem:3}
For small enough $\xi$, $Q_0$ almost surely, under $P^\o_0$, the
family of processes $(X^{\xi,\eps}(t)=\eps X^\xi(\frac
t{\eps^2}),t\in\R_+)$ is tight in the Skorokhod topology (as $\eps$ goes to $0$).
\end{lm}

{\it Proof}: we shall prove that, for any $T>0$, for any $\eta>0$ and for small enough $\theta_0$ then  
\beqn \label{eq:ti10} \limsup_\eps 
\sup_{v\leq T} P^\o_0(\sup_{u\leq T\,;\, v\leq u\leq v+\theta_0} \vert
X^{\xi,\eps}(u)-X^{\xi,\eps}(v)\vert >\eta) \leq
27 (c_e)^3 (\frac{\sqrt
{\theta_0}}\eta)^3\,. \eeqn 
Indeed inequality (\ref{eq:ti10}) implies that 
 \beqn \label{eq:ti11} 
 \limsup_{\theta_0} \frac 1 {\theta_0} \limsup_\eps
\sup_{v\leq T}P^\o_0(\sup_{u\leq T\,;\, v\leq u\leq v+\theta_0} \vert
X^{\xi,\eps}(u)-X^{\xi,\eps}(v)\vert >\eta)=0\,.\eeqn 
According to Theorem 8.3 in Billingsley's book \cite{kn:Bill},  this 
last inequality is sufficient to ensure the tightness. 

We use Lemma \ref{lem:2} with $\theta=1$ to check
that 
\beqnn P^\o_0(\sup_{t\leq T}\vert
X^{\xi,\eps}(t)\vert\geq K) =P^\o_0(\tau(0,\frac K\eps)\leq\frac
T{\eps^2})\leq c_e\frac{\sqrt T}K\,. \eeqnn 
(We could use Lemma
\ref{lem:2} since $\frac K\eps\leq \frac T{\eps^2}$ for small
$\eps$.) 

Next choose $\eta>0$ and use Lemma \ref{lem:2'} with $\theta=3$ and the Markov property to get
that \beqnn P^\o_0(\sup_{v\leq u\leq T\,;\, u-v\leq\theta_0} \vert
X^{\xi,\eps}(u)-X^{\xi,\eps}(v)\vert >\eta) \leq
P^\o_0(&&\sup_{t\leq T}\vert X^{\xi,\eps}(t)\vert\geq K) \\&&+
\sup_{y\,;\, \vert y\vert\leq K/\eps}
P^\o_y(\tau(y,\frac\eta\eps)\leq \frac{\theta_0}{\eps^2})\,.\eeqnn
If we choose $K$ of order $1/\eps$ and pass to the limit as $\eps$ tends to $0$, 
then, due to the previous inequality, the contribution of the first term vanishes. 
As for the second term, by Lemma \ref{lem:2'}, it is bounded by 
$27 (c_e)^3 (\frac{\sqrt
{\theta_0}}\eta)^3$. Note that we could use 
Lemma
\ref{lem:2'} since $\frac K\eps\leq (\frac {\theta_0}{\eps^2})^3$ 
and $\frac \eta\eps\leq (\frac {\theta_0}{\eps^2})^3$  for small
$\eps$. Thus the proof of (\ref{eq:ti10}) is complete. \qed

\vskip .5cm {\bf Step 4: Poincar\'e inequalities and end of the
proof of Theorem \ref{theorem'}}

Applied to a centered function $f$, 
Nash inequality (\ref{eq:nash}) 
reads: 
\beqnn 
 \Vert f\Vert_2^{2+\frac{4}{\eps(n)}}
  \leq
\beta\,
  n^{2(1-\frac{d}{\eps(n)})}\,\tEEon(f,\,f)\,\|f\|_1^{4/\eps(n)}\,.
\eeqnn

Holder's inequality implies that 
$$\Vert f\Vert_1\leq \Vert f\Vert_2 (2n+1)^{d/2}$$ since $\# C^n(\o)\leq (2n+1)^d$. 
We deduce that any centered function on $\CC^n(\o)$ satisfies 
\beqnn \Vert f\Vert_2^2 \leq \beta  n^2\,\tEEon(f,\,f)\,, \eeqnn
for some constant $\beta$. Equivalently,  any (not necessarily centered) 
function on $\CC^n(\o)$ satisfies 
\beqnn \Var (f)\leq \beta  n^2\,\tEEon(f,\,f)\,. \eeqnn

Thus we have proved the following Poincar\'e inequality on $C^n(\o)$: there is a constant $\beta$
such that, $Q_0^\xi$.a.s. for large enough $n$, for any function
$f:\CC^n(\o)\rightarrow\R$ then \beqn\label{eq:poinc1} \sum_{x\in
\CC^n(\o)}f(x)^2 \leq \beta  n^2\, \sum_{x\sim y\in
\CC^n(\o)}\oxi(x,y) (f(x)-f(y))^2 \eeqn

Our second Poincar\'e inequality is derived from \cite{kn:Ba}, see
Definition 1.7, Theorem 2.18, Lemma 2.13 part a) and Proposition 2.17 part b): there
exist constants $M<1$ and $\beta$ such that $Q_0^\xi$.a.s. for any
$\delta>0$, for large enough $n$, for any $z \in \Z^d$ s.t. $\vert
z\vert\le n$ and for any function $f:\Z^d\rightarrow\R$ then
\beqn\label{eq:poinc2} \sum_{x\in \CC^\xi(\o)\cap (z+[-M\delta
n,M\delta n]^d)}
                        f(x)^2
\leq \beta \delta^2 n^2\, \sum_{x\sim y\in \CC^\xi(\o)\cap
(z+[-\delta n,\delta n]^d)}
                        \oxi(x,y) (f(x)-f(y))^2
\eeqn In \cite{kn:Ba}, inequality (\ref{eq:poinc2}) is in fact
proved for the Dirichlet form $\EE^{\alpha',n}$ but the comparison
inequality (\ref{eq:compDF}) implies that it also holds for the
Dirichlet form $\tEEon$.

One can now conclude the proof of the Theorem following the argument
in \cite{kn:MP} line by line starting from paragraph 2.2. \qed


\vskip 1cm
\section{Percolation results}
\setcounter{equation}{0}
\label{sec:perco}

\subsection{Prerequisites on site percolation}

We shall use some properties of site percolation that we state below.

By site percolation of parameter $r$ on $\Z^d$, we mean the product
Bernoulli measure of parameter $r$ on the set of applications
$\zeta:\Z^d\rightarrow\{0,1\}$. We identify any such application
with the sub-graph of the grid whose vertices are the points
$x\in\Z^d$ such that $\zeta(x)=1$ and equipped with the edges of the
grid linking two points $x,y$ such that $\zeta(x)=\zeta(y)=1$.

Let $l>1$. Call a sub-set of $\Z^d$ {\it $l$-connected} if it is
connected for the graph structure defined by: two points are
neighbors when the Euclidean distance between them is less than $l$.

We recall
our notation  $\vert x-y\vert$ for the Euclidean distance between $x$ and $y$.

A {\it path} is a sequence of vertices of $\Z^d$ such that 
two successive vertices in $\pi$ are neighbors. 
We mostly consider injective paths. 
With some abuse of vocabulary, a sequence of vertices of $\Z^d$ 
in which two successive vertices are at distance not more 
than $l$ will be called a {\it $l$-nearest-neighbor path}. 
Let $\pi=(x_0,...,x_k)$ be a sequence of vertices. We define its length 
$$\vert \pi\vert=\sum_{j=1}^k \vert x_{j-1} -x_j\vert\,,$$ 
and its cardinality $\#\pi=\#\{x_0,...,x_k\}$. ($\#\pi=k+1$ for an injective path.) 
When convenient, we identify an injective path with a set (its range).

\begin{lm}\label{lem:site}
Let $l>1$.
There exists $p_1>0$ such that for $r<p_1$,
almost any realization of site percolation of parameter $r$ has only
finite $l$-connected components and,
for large enough $n$, any $l$-connected component that
intersects the box $[-n,n]^d$ has volume smaller than $(\log n)^{6/5}$.
\end{lm}

{\it Proof}: the number of $l$-connected sets that contain a fixed
vertex and of volume $m$ is smaller than $e^{a(l)m}$ for some
constant $a(l)$, see \cite{kn:G}. Thus the number of $l$-connected
sets of volume $m$ that intersect the box $[-n,n]^d$ is smaller than
$(2n+1)^d e^{a(l)m}$. But the probability that a given set of volume
$m$ contains only opened sites is $r^m\leq p_1^m$. We now choose
$p_1$ small enough so that $\sum_n\sum_{m\geq (\log n)^{6/5}}
(2n+1)^d e^{a(l)m} p_1^m<\infty$ and the Borel-Cantelli lemma yields
the conclusion of Lemma \ref{lem:site}. \qed

As in the case of bond percolation discussed in the introduction, it
is well known that for $r$ larger than some critical value then
almost any realization of site percolation of parameter $r$ has a
unique infinite connected component - the {\it infinite cluster} -
that we will denote with $\CC$.

\begin{lm}\label{lem:site'}
There exists $p_2<1$ such that for $r>p_2$, for
almost any realization of site percolation of parameter $r$  and
for large enough $n$, any connected component of the complement of the
infinite cluster $\CC$ that
intersects the box $[-n,n]^d$ has volume smaller than $(\log n)^{5/2}$.
\end{lm}

{\it Proof}: let $\zeta$ be a typical realization of site percolation of parameter
$r$.
We assume that $r$ is above the critical value so that there is a unique infinite
cluster, $\CC$.
We also assume that $1-r<p_1$ where $p_1$ is the value provided by Lemma
\ref{lem:site} for $l=d$.

Let $A$ be a connected component of the complement of $\CC$. Define
the {\it interior boundary of $A$}: $\partial_{int}A=\{x\in A\,;\,
\exists y\, s.t.\, (x,y)\in \E_d\,\hbox{and}\, y\notin A\}$. 
It is known  
that $\partial_{int}A$ is $d$-connected, 
 see \cite{kn:DP}, Lemma 2.1. 
By construction any $x\in
\partial_{int}A$ satisfies $\zeta(x)=0$. Since the application
$x\rightarrow 1-\zeta(x)$ is a typical realization of site
percolation of parameter $1-r$ and $1-r<p_1$, as an application of
Lemma \ref{lem:site} we get that $\partial_{int}A$ is finite.
Because we already know that the complement of $A$ is infinite
(since it contains $\CC$), it implies that $A$ itself is finite.

We now assume that $A$ intersects the box $[-n,n]^d$. Choose $n$
large enough so that $\CC\cap [-n,n]^d\not=\emptyset$ so that
$[-n,n]^d$ is not a sub-set of $A$. Then it must be that
$\partial_{int}A$ intersects $[-n,n]^d$. Applying Lemma
\ref{lem:site} again, we get that, for large $n$, the volume of
$\partial_{int}A$ is smaller than $(\log n)^{6/5}$. The classical
isoperimetric inequality in $\Z^d$ implies that, for any finite
connected set $B$, one has $(\#\partial_{int} B)^{d/(d-1)}\geq
{\cal I} \# B$ for some constant $\cal I$. Therefore $\# A\leq
{\cal I}^{-1} ( \log n)^{6d/5(d-1)}$. Since $6d/5(d-1)<5/2$, the
proof is complete. \qed

\begin{lm}\label{lem:site''}
There exists $p_3<1$ and a constant $c_3$ such that for $r>p_3$, for
almost any realization of site percolation of parameter $r$  and
for large enough $n$, for any two points $x,y$ in the box $[-n,n]^d$ 
such that $\vert x-y\vert\geq (\log n)^{3/2}$ we have\\ 
(i) for any injective $d$-nearest-neighbor path $\pi$ from $x$ to $y$ then  
\beqnn 
\#\{z\in\pi\,;\,\zeta(z)=1\}\geq c_3\vert x-y\vert\,. 
\eeqnn
(ii) for any injective ($1$-nearest-neighbor) path $\pi$ from $x$ to $y$ then  
\beqnn 
\#(\CC\cap\pi)\geq c_3\vert x-y\vert\,. 
\eeqnn
\end{lm}

{\it Proof}: we assume that $r$ is close enough to $1$ so that there is a unique 
infinite cluster $\CC$. We also assume that $1-r<p_1$, 
where $p_1$ is the constant appearing in Lemma \ref{lem:site} for $l=1$.  
Then the complement 
of $\CC$ only has finite connected components. 

Part (i) of the Lemma is proved by a classical Borel-Cantelli argument 
based on the following simple observations: the number of 
injective $d$-nearest-neighbor paths $\pi$ from $x$ of length $L$ is bounded by $(c_d)^L$ for 
some constant $c_d$ that depends on the dimension $d$ only; the probability that 
a given set of cardinality $L$ contains less than $dc_3L$ sites where $\zeta=1$ is bounded 
by $exp(\lambda dc_3 L)(re^{-\lambda}+1-r)^L$ for all $\lambda>0$. 
We choose $c_3<\frac 1d$ 
and $\lambda$ such that $c_d e^{-(1-dc_3)\lambda}<1$ and $p_3$ such that 
$\gamma=c_d e^{\lambda dc_3} (p_3e^{-\lambda}+1-p_3)<1$. 
Let now $x$ and $y$ be as in the Lemma. Note that any 
injective $d$-nearest-neighbor path $\pi$ from $x$ to $y$ 
satisfies $\#\pi\geq \frac 1d \vert x-y\vert\geq \frac 1d (\log n)^{3/2}$. 
Therefore the probability that there is an injective $d$-nearest-neighbor 
path $\pi$ from $x$ to $y$ such that $\#\{z\in\pi\,;\,\zeta(z)=1\}<c_3\vert x-y\vert$ 
is smaller than 
$\sum_{L\geq \frac 1d (\log n)^{3/2}} \gamma^L$ and the probability 
that (i) fails for some $x$ and $y$ is smaller than 
$(2n+1)^{2d}\sum_{L\geq \frac 1d (\log n)^{3/2}} \gamma^L$. 
Since $\sum_n (2n+1)^{2d} \sum_{L\geq \frac 1d (\log n)^{3/2}} \gamma^L<\infty$, 
the Borel-Cantelli lemma then yields that, for large enough $n$, part (i) of Lemma 
\ref{lem:site''} holds. 

We prove part (ii) by reducing it to an application of part (i). 
Assume that, for some points $x$ and $y$ as in the Lemma, there exists 
an injective 
nearest-neighbor path $\pi$ from $x$ to $y$ such that   
$\#(\CC\cap\pi)< c_3\vert x-y\vert$. 
We first modify the path $\pi$ into a $d$-nearest-neighbor path from $x$ to $y$, say  
$\pi'$, in the following way: the parts of $\pi$ that lie in $\CC$ remain unchanged but 
the parts of $\pi$ that visit the complement of $\CC$ are  modified so that they only 
visit points where $\zeta=0$. Such a modified path $\pi'$ exists because 
the interior boundary of a connected component of the complement of $\CC$ 
is $d$ connected (as we already mentioned in the proof of Lemma \ref{lem:site'}) 
and only contains points where $\zeta=0$. 

Observe that $\CC\cap\pi'=\CC\cap\pi$ and that $\CC\cap\pi'=\{z\in\pi'\,;\,\zeta(z)=1\}$ 
so that 
\beqnn \#\{z\in\pi'\,;\,\zeta(z)=1\}< c_3\vert x-y\vert\,.
\eeqnn 
Next turn $\pi'$ into an injective $d$-nearest-neighbor path, say $\pi''$, by suppressing loops 
in $\pi'$. Clearly $\{z\in\pi''\,;\,\zeta(z)=1\}\subset\{z\in\pi'\,;\,\zeta(z)=1\}$ 
and therefore 
\beqnn \#\{z\in\pi''\,;\,\zeta(z)=1\}< c_3\vert x-y\vert\,,
\eeqnn 
a contradiction with part (i) of the Lemma. \qed

\subsection{Proof of Lemma \ref{lem:holes}}

Lemma \ref{lem:holes} only deals with the geometry of percolation
clusters, with no reference to random walks. We will restate it as a
percolation lemma at the cost of changing a little our notation. In
order to make a distinction with a typical realization of an
environment for which we used the notation $\omega$, we will use the
letters $\alpha$ or $\alpha'$ to denote typical realizations of a
percolation graphs. Thus one switches from the notation of the
following proof back to the notation of part \ref{sec:deduce} using
the following dictionary: \beqnn
\alpha(e)=\1_{\o(e)>0}\,&;&\, \alpha'(e)=\1_{\o(e)\geq\xi}\\
q=Q(\o(e)>0)\,&;&\, p=Q(\o(e)\geq\xi\,\vert\, \o(e)>0)\,.
\eeqnn
This way taking $\xi$ close to $0$ is equivalent to taking $p$ close to $1$.

We very much rely on renormalization technics, see Proposition 2.1.
in \cite{kn:AP}.

As in the introduction, we identify a sub-graph of $\Z^d$ with an
application $\alpha:\E_d\rightarrow\{0,1\}$, writing $\alpha(x,y)=1$
if the edge $(x,y)$ is present in $\alpha$ and $\alpha(x,y)=0$
otherwise. Thus $\AAA=\{0,1\}^{\E_d}$ is identified with the set of
sub-graphs of $\Z^d$. Edges pertaining to $\a$ are then called {\it
open}. Connected components of such a sub-graph will be called {\it
clusters}.

Define now $Q$ to be the probability measure on $\{0,1\}^{\E_d}$ under which
the random
variables $(\alpha(e),\,e \in \E_d)$ are Bernoulli$(q)$ independent variables
with
\beqnn q>p_c.
\eeqnn
Then, $Q$ almost surely, the graph $\a$ has a unique infinite cluster denoted with
$\CC(\a)$.

For a typical realization of the percolation graph under $Q$, say $\alpha$, let
$Q^\alpha$ be the law of bond percolation on $\CC(\alpha)$ with parameter $p$.
We shall denote $\alpha'$ a typical realization under $Q^\alpha$ i.e. $\alpha'$
is a random subgraph of $\CC(\alpha)$ obtained by keeping (resp. deleting) edges
with probability $p$ independently of each other.
We always assume that $p$ is close enough to $1$ so that $Q^\alpha$ almost surely
there is a unique infinite cluster in $\alpha'$ that we denote $\CC^\alpha(\alpha')$.
By construction $\CC^\alpha(\alpha')\subset \CC(\alpha)$.
Connected components of the complement of $\CC^\alpha(\alpha')$ in $\CC(\alpha)$ are
called {\it holes}.


We now restate Lemma \ref{lem:holes}:\\
{\it there exists $p_0<1$ such that for $p>p_0$, for $Q$ almost any $\alpha$, for
$Q^\alpha$ almost any $\alpha'$, for large enough $n$, then any hole
intersecting the box $[-n,n]^d$ has volume smaller than $(\log n)^a$.}

{\it Renormalization}: let $\alpha$ be a typical realization of percolation under $Q$.

Let $N$ be an integer. We chop $\Z^d$ in a disjoint union of boxes
of side length $2N+1$. Say $\Z^d=\cup_{{\mathbf i}\in\Z^d}B_{\mathbf
i}$, where $B_{\mathbf i}$ is the box of center $(2N+1){\mathbf i}$.
Following \cite{kn:AP}, let $B'_{\mathbf i}$ be the box of center
$(2N+1){\mathbf i}$ and side length $\frac 5 2 N +1$. From now on,
the word {\it box} will mean one of the boxes $B_{\mathbf i},
{\mathbf i}\in\Z^d$.

We say that a box $B_{\mathbf i}$ is {\it white} if
$B_{\mathbf i}$ contains at least one edge from $\alpha$
and the event $R_{\mathbf i}^{(N)}$ in equation (2.9) of
\cite{kn:AP} is satisfied. Otherwise, $B_{\mathbf i}$ is a {\it black} box.
We recall that the event $R_{\mathbf i}^{(N)}$ is defined by:
there is a unique cluster of $\alpha$ in $B'_{\mathbf i}$, say $K_{\mathbf i}$;
all open paths
contained in $B'_{\mathbf i}$ and of radius larger than $\frac 1 {10} N$
intersect $K_{\mathbf i}$ within $B'_{\mathbf i}$; $K_{\mathbf i}$ is crossing for
each subbox
$B\subset B'_{\mathbf i}$ of side larger than $\frac 1 {10} N$.
See \cite{kn:AP} for details.
We call $K_{\mathbf i}$ the {\it crossing cluster }
of $\alpha$ in the box $B_{\mathbf i}$.
Note the following consequences of this definition.

(Fact i) If $x$ and $y$ belong to the same white box $B_{\mathbf i}$ and
both $x$ and $y$ belong to the infinite cluster of $\alpha$, then
there is a path in $\CC(\alpha)$ connecting $x$ and $y$ within
$B'_{\mathbf i}$.

(Fact ii) Choose two neighboring indices $\mathbf i$ and $\mathbf j$
with $\vert {\mathbf i}-{\mathbf j}\vert=1$ and such that both boxes
$B_{\mathbf i}$ and $B_{\mathbf j}$ are white. As before, let
$K_{\mathbf i}$ and $K_{\mathbf j}$ be the crossing clusters in
$B_{\mathbf i}$ and $B_{\mathbf j}$ respectively. Let $x\in
K_{\mathbf i}$ and $y\in K_{\mathbf j}$. Then there exists a path in
$\alpha$ connecting $x$ and $y$ within $B'_{\mathbf i}\cup
B'_{\mathbf j}$.

We call {\it renormalized} process the random subsets of $\Z^d$ obtained
by taking the image of the initial percolation model by the application $\phi_N$,
see equation (2.11) in \cite{kn:AP}. A site $\mathbf i \in\Z^d$ is thus declared
{\it white} if the box $B_{\mathbf i}$ is white.

Let $\mathbf Q$ be the law of the renormalized process. The comparison result of
Proposition 2.1 in \cite{kn:AP} states that $\mathbf Q$ stochastically dominates
the law of site percolation with parameter $p(N)$ with
$p(N)\rightarrow 1$ as $N$ tends to $\infty$.

We now introduce the extra percolation $Q^\alpha$. Let us call {\it grey} a white
box $B_{\mathbf i}$ that contains
an edge $e\in \CC(\alpha)$ such that $\alpha'(e)=0$. We call
{\it pure white} white boxes that are not grey.

Let $\mathbf Q'$ be the law on subsets of the renormalized grid
obtained by keeping pure white boxes, and deleting both black and
grey boxes. We claim that $\mathbf Q'$ dominates the law of site
percolation with parameter $p'(N)=p(N) p^{\,e_N(d)}$ where $e_N(d)$
is the number of edges in a box of side length $2N+1$. (Remember
that $p$ is the parameter of $Q^\alpha$.) This claim is a
consequence of the three following facts. We already indicated that
$\mathbf Q$ stochastically dominates the law of site percolation
with parameter $p(N)$. The conditional probability that a box
$B_{\mathbf i}$ is pure white given it is white is larger or equal
than $p^{\,e_N(d)}$. Besides, still under the condition that
$B_{\mathbf i}$ is white, the event `$B_{\mathbf i}$ is pure white'
is independent of the colors of the other boxes.

We further call {\it immaculate} a pure white box $B_{\mathbf i}$
such that any box $B_{\mathbf j}$ intersecting  $B'_{\mathbf i}$ is
also pure white. Call $\mathbf Q''$ the law on subsets of the
renormalized grid obtained by keeping only immaculate boxes. Since
the event `$B_{\mathbf i}$ is immaculate' is an increasing function
with respect to the percolation process of pure white boxes, we get
that $\mathbf Q''$ stochastically dominates the law of site
percolation with parameter $p''(N)=p'(N)^{3^d}$.

{\it End of the proof of Lemma \ref{lem:holes}}: choose $p_0$ and $N$
such that $p''(N)$ is close enough to $1$
so that, $\mathbf Q''$ almost surely,
there is an infinite cluster of immaculate boxes that we call $\C$.

For $\mathbf i\in\C$, let $K_{\mathbf i}$ be the crossing cluster in
the box $B_{\mathbf i}$ and let $K=\cup_{\mathbf i\in\C}K_{\mathbf i}$.
Then $K$ is connected (This follows from the definition of white
boxes, see (Fact i) and (Fact ii) above.) and infinite (Because $\C$
is infinite.). Thus we have $K\subset \CC^\alpha(\alpha')$.

Let $A$ be a hole and let $\mathbf A$ be the set of indices $\mathbf i$ such that
$B_{\mathbf i}$
intersects $A$. Observe that $\mathbf A$ is connected.
We claim that $$\mathbf A\cap\C=\emptyset\,.$$

Indeed, assume there exists $x\in B_{\mathbf i}$ such that $\mathbf
i\in \C$ and $x\in A$. By definition $A$ is a subset of
$\CC(\alpha)$ and therefore $x\in\CC(\alpha)$. Let $y\in K_{\mathbf
i}$, $y\not=x$. As we already noted $y\in\CC^\alpha(\alpha')$. Since
$x\in \CC(\alpha)$ and  $y\in\CC(\alpha)$ there is a path, $\pi$,
connecting $x$ and $y$ within $B'_{\mathbf i}$, see (Fact i) above.
But $B_{\mathbf i}$ is immaculate and therefore $B'_{\mathbf i}$
only contains edges $e$ with $\alpha'(e)=1$. Therefore all edges
along the path $\pi$ belong to $\alpha'$ which imply that 
$x\in\CC^\alpha(\alpha')$. This is in contradiction
with the assumptions that $x\in A$. 
We have proved that $\mathbf A\cap\C=\emptyset$.

To conclude the proof of Lemma \ref{lem:holes}, it only remains to choose $p_0$ and
$N$
such that $p''(N)\geq p_2$ and apply Lemma \ref{lem:site'}. We deduce that the
volume of
$\mathbf A$ is bounded by $(\log n)^{5/2}$ and therefore the volume of $A$ is smaller
than $(2N+1)^d (\log n)^{5/2}$.
\qed

\subsection{Deviation of the chemical distance} \label{sec:deviation}

We use the same notation as in the preceeding section. For given
realizations of the percolations $\alpha$ and $\alpha'$, we define  
the corresponding {\it chemical distance} $d^\alpha_{\alpha'}$ on 
$\CC^\alpha(\alpha')$: two points $x\not=y$ in  $\CC^\alpha(\alpha')$ 
satisfy $d^\alpha_{\alpha'}(x,y)=1$ if and only if one (at least) of the following 
two conditions is satisfied: either $x$ and $y$ are neighbors in $\Z^d$ and 
$\alpha'(x,y)=1$ or both $x$ and $y$ are at the boundary of a hole $h$ 
i.e. there is a hole $h$ and $x',y'\in h$ such that $x'$ is a neighbor 
of $x$ and $y'$ is a neighbor of $y$. 
In general, $d^\alpha_{\alpha'}(x,y)$ is defined as the smaller integer $k$ 
such that there exists a sequence of points $x_0,...,x_k$ in $\CC^\alpha(\alpha')$ 
with $x_0=x$, $x_k=y$ and 
such that $d^\alpha_{\alpha'}(x_j,x_{j+1})=1$ for all $j$.

\begin{lm}\label{lem:deviation}
There exists $p_4<1$ such that for $p>p_4$, there exist constants
$c^+$ and $c^-$ such that for $Q$ almost any $\alpha$, for
$Q^\alpha$ almost any $\alpha'$, for large enough $n$, then
\beqn\label{eq:deviation} c^- \vert x-y\vert\leq
d^\alpha_{\alpha'}(x,y) \leq c^+ \vert x-y\vert\,,\eeqn for any
$x,y\in\CC^\alpha(\alpha')$ such that $x\in[-n,n]^d$ and $\vert
x-y\vert\geq (\log n)^2$.
\end{lm}

{\it Proof}: let $d^\alpha(x,y)$ be the chemical distance between
$x$ and $y$ within $\CC(\alpha)$ i.e. $d^\alpha(x,y)$ is the minimal
length of a path from $x$ to $y$, say $\pi$, such that any edge
$e\in\pi$ satisfies $\alpha(e)=1$. \\
Applying Theorem 1.1 in \cite{kn:AP} together with the
Borel-Cantelli Lemma, we deduce that there exists a constant $c^+$
such that $d^\alpha(x,y) \leq c^+ \vert x-y\vert$ for any
$x,y\in\CC(\alpha)$ such that $x\in[-n,n]^d$ and $\vert x-y\vert\geq
(\log n)^2$.  Since $d^\alpha_{\alpha'}(x,y) \leq d^\alpha(x,y)$, it gives
the upper bound in (\ref{eq:deviation}).

We now give a proof of the lower bound. 
As for Lemma \ref{lem:holes}, we use a renormalization argument. 
The notation used below is borrowed from the proof of Lemma \ref{lem:holes} 
except that the role of $p_0$ is now played by $p_4$. . 

We wish to be able to apply Lemma \ref{lem:site''} (ii) to the renormalized site percolation model 
with law $\mathbf Q''$ (i.e. the percolation model of immaculate boxes): 
therefore we choose $p_4$ and $N$ such that $p''(N)\geq p_3$ and observe that the event considered in 
Lemma \ref{lem:site''} (ii) is increasing. 

Consider two points $x$ and $y$ as in Lemma \ref{lem:deviation} 
and let $\pi$ be an injective path from $x$ to $y$ within $\CC(\alpha)$. We shall prove that 
\beqn\label{eq:5.3.1} 
\#\EE_\pi\geq c_5\vert x-y\vert\,,\eeqn 
where $\EE_\pi=\{z,z'\in\pi\cap\CC^\alpha(\alpha')\,;\, \alpha'(z,z')=1\}$. 
By construction of the chemical distance $d^\alpha_{\alpha'}$, 
(\ref{eq:5.3.1}) implies the lower bound in (\ref{eq:deviation}) 
with $c^-=c_5$. 

Let $\Pi'$ be the sequence of the indices of the boxes $B_{\mathbf i}$ that $\pi$ 
intersects. At the level of the renormalized grid, $\Pi'$ is a nearest-neighbor 
path from $\mathbf i_0$ to $\mathbf i_k$ with $x\in B_{\mathbf i_0}$ 
and $y\in B_{\mathbf i_k}$. 
Let $\Pi=({\mathbf i}_0,...,{\mathbf i}_k)$ be the injective path 
obtained by suppressing loops in $\Pi'$.  
We may, and will, assume that $n$ is large enough so that 
$i_0\not= i_k$ so that $\vert \mathbf i_0-\mathbf i_k \vert$ and $\vert x-y\vert$ are comparable. 
Applying Lemma \ref{lem:site''} (ii) to $\mathbf Q''$, we get that 
\beqn\label{eq:5.3.2}\# (\C\cap\Pi)\geq c_3\vert \mathbf i_0-\mathbf i_k \vert
\geq c'_3\vert x-y\vert\,,\eeqn 
for some constant $c'_3$.

Let ${\mathbf i}\in \C\cap\Pi$ and choose $z\in B_{\mathbf i}\cap\pi$. 
Since the path $\pi$ is not entirely contained in one box, 
it must be that $\pi$ connects $z$ to some point 
$z'\notin B_{\mathbf i}$. Since $z'\in\pi$, we also have $z'\in \CC(\alpha)$. 
By definition of a white box, it 
implies that $z\in K_{\mathbf i}$. 
Since ${\mathbf i}\in \C$, it implies that actually $z\in K$ and therefore 
$z\in \CC^\alpha(\alpha')$. As a matter of fact, since the box $B_{\mathbf i}$ is pure  
white, we must have $\alpha'=1$ on all the edges of $\pi$ from $z$ to $z'$. 
In particular  $z$ has a neighbor in $\CC(\alpha)$, say $z''$, 
such that $\alpha'(z,z'')=1$. Therefore $(z,z'')\in\EE_\pi$. 
We conclude that any indice in $\C\cap\Pi$ 
gives a contribution of at least $1$ to $\#\EE_\pi$. 
Therefore (\ref{eq:5.3.2}) implies that  $$ \#\EE_\pi\geq c'_3\vert x-y\vert\,.$$ 
\qed



\begin{thebibliography}{xxxxxx 89}

\bibitem{kn:AP} Antal P.,Pisztora A.~(1996)\\
                On the chemical distance for supercritical Bernouilli percolation\\
                {\em Ann.~Probab.}~{\bf 24}, 1036-1048.

\bibitem{kn:Ba} Barlow, M.T.~(2004)\\
                Random walks on supercritical percolation clusters\\
                {\em Ann.~Probab.}~{\bf 32}, 3024-3084.
                
 \bibitem{kn:Bass} Bass, R.F.~(2002)\\
                On Aronson's upper bounds for heat kernels.\\
                {\em Bull.~London~Math.~Soc.}~{\bf 34}, 415-419.
                               

\bibitem{kn:BB}   Berger, N., Biskup, M.~(2007)\\
                Quenched invariance principle for simple random walk
                on percolation clusters.  \\ 
                {\em  Prob.~Th.~Rel.~Fields}~{\bf 137}, 83-120.  
                
\bibitem{kn:BBHK}   Berger, N., Biskup, M., Hoffman, C., Kozma, G.~(2006)\\
                Anomalous heat kernel decay for random walk among bounded random conductances.  
                To appear in  {\em  Ann.~Inst.~Henri~Poincar\'e}.     
                
\bibitem{kn:Bill}    Billingsley, P.~(1968)  \\                 
                {\em The convergence of probability measures.}\\ 
                John Wiley, New York.
                
\bibitem{kn:BiPres}  Biskup, M., Prescott, T.M.~(2007)\\ 
		Functional CLT for random walk among bounded random conductances.\\ 
		Preprint 2007.  


               
		
\bibitem{kn:DFGW} De Masi, A.,  Ferrari, P., Goldstein, S.,  Wick, W.D.~(1989)\\
                  An invariance principle for reversible Markov
                  processes. Applications to random motions in random environments\\
                 {\em Journ.~Stat.~Phys.}~{\bf 55}~(3/4),  787-855.
                 
\bibitem{kn:DP} Deuschel, J-D., Pisztora A.~(1996)\\ 
		Surface order deviations for high density percolation.\\ 
		{\em Prob.~Th.~Rel.~Fields}~{\bf 104}, 467-482.
		  

\bibitem{kn:EK} Ethier, S.N., Kurtz, T.G.~(1986)\\
                {\em Markov processes}\\
                John Wiley, New York.

\bibitem{kn:FM} Fontes, L.R.G., Mathieu, P.~(2006)\\
                On symmetric random walks with random conductances on $\Z^d$\\
                {\em Prob.~Th.~Rel.~Fields}~{\bf 134}, 565-602.


\bibitem{kn:G}  Grimmett, G.~(1999)\\
                {\em Percolation}\\
                Springer-Verlag, Berlin (Second edition).

 \bibitem{kn:JS}  Jacod, J., Shiryaev, A.N.~(1987)\\
                {\em Limit theorems for stochastic processes}\\
                Springer-Verlag, Berlin.

\bibitem{kn:Ko} Kozlov, S.M.~(1985)\\
                The method of averaging and walks in inhomogeneous environments\\
                {\em Russian~Math.~Surveys}~{\bf 40}~(2), 73-145.


\bibitem{kn:MP} Mathieu, P.,  Piatnitski, A.L.~(2007)\\
               Quenched invariance principles for random walks on percolation clusters.\\  
                {\em Proc.~R.~Soc.~A}~{\bf 463}, 2287-2307. 


\bibitem{kn:MR} Mathieu, P., Remy, E.~(2004)\\
                Isoperimetry and heat kernel decay on percolations clusters\\
                {\em Ann.~Probab.}~{\bf 32}, 100-128.
                
\bibitem{kn:Nash} Nash, J.~(1958)\\
                Continuity of solutions of parabolic and elliptic equations\\
                {\em Amer.~J.~Math.}~{\bf 80}, 931-954.
                                

\bibitem{kn:SS} Sidoravicius, V., Sznitman, A-S.~(2004)\\
                Quenched invariance principles for walks on clusters of percolation
                or among random conductances\\
                {\em Prob.~Th.~Rel.~Fields}~{\bf 129}, 219-244.








\end{thebibliography}
\end{document}